\documentclass{article}
\usepackage{spconf,amsmath,graphicx,hyperref,amsfonts, amssymb,bm,bbm}
\usepackage{color,soul}
\usepackage{algorithm,algpseudocode}
\usepackage[square,numbers]{natbib}

\newcommand*\w{ {\boldsymbol{\scriptstyle \mathcal{W} }} }
\newcommand*\A{\mathcal{A}}
\newcommand*\B{\mathcal{B}}
\newcommand*\transp{^\mathsf{T}}
\newcommand*\J{\mathcal{J}}
\newcommand*\werr{\widetilde{\w}}
\newcommand*\gradcal{\nabla\J}
\newcommand*\gradappcal{\widehat{\gradcal}}

\newcommand\norm[1]{%
\left\lVert#1\right\rVert
}
\allowdisplaybreaks

\newtheorem{assumption}{Assumption}

\title{Graph-Aware Learning Rates for Decentralized Optimization}

\name{Aaron Fainman and Stefan Vlaski\thanks{ This work was supported in part by EPSRC Grants EP/X04047X/1 and EP/Y037243/1. Emails: \{aaron.fainman22,s.vlaski\}@imperial.ac.uk.}}
\address{Department of Electrical and Electronic Engineering, Imperial College London}

\begin{document}
%
\maketitle
\begin{abstract}
We propose an adaptive step-size rule for decentralized optimization. Choosing a step-size that balances convergence and stability is challenging. This is amplified in the decentralized setting as agents observe only local (possibly stochastic) gradients and global information (like smoothness) is unavailable. We derive a step-size rule from first principles. The resulting formulation reduces to the well-known Polyak's rule in the single-agent setting, and is suitable for use with stochastic gradients. The method is parameter free, apart from requiring the optimal objective value, which is readily available in many applications. Numerical simulations demonstrate that the performance is comparable to the optimally fine-tuned step-size.
\end{abstract}


\section{Introduction}
Consider a network of $K$ agents aiming to solve a consensus optimization problem of the form:
\begin{align}\label{eq:consensus-opt}
    w^o=\underset{w\in\mathbb{R}^M}{\arg\min}\left\{J(w)\triangleq\frac{1}{K}\sum_{k=1}^K J_k(w)\right\}
\end{align}
where $J_k(w)\triangleq\mathbb{E} Q_k(w;\boldsymbol{x}_k)$ denotes the local objective function of agent $k$. Here, \( Q_k(\cdot; \cdot) \) denotes the loss function, \( w \in \mathbb{R}^M \) denotes the parameters to be chosen, and \( \boldsymbol{x}_k \) denotes the local data available to agent \( k \).

Many algorithms exist to solve problem~\eqref{eq:consensus-opt} (see, for example~\cite{sayed14,vlaski23}), \textcolor{black}{ including gradient-tracking~\cite{dilorenzo16,augdgm} and primal-dual methods~\cite{extra, exact_diff}. Here we focus on ATC-Diffusion~\cite{lopes08,chen15}}:
\begin{subequations}\label{eq:atc}
\begin{align}
    \boldsymbol{g}_{k,i} &= \boldsymbol{w}_{k,i-1} - \boldsymbol{\mu}_{k,i-1}\widehat{\nabla J}_k(\boldsymbol{w}_{k,i-1})\label{eq:atc-opt}\\ 
    \boldsymbol{w}_{k,i} &= \sum_{\ell\in\mathcal{N}_k} a_{k\ell}\boldsymbol{g}_{\ell,i} \label{eq:atc_avg}
\end{align}
\end{subequations}
Each agent updates its parameter vector $\boldsymbol{w}_{k,i}$ using a stochastic gradient estimate in~$\eqref{eq:atc-opt}$, then averages with its neighbors~$\mathcal{N}_k$ in~$\eqref{eq:atc_avg}$. A common choice for the gradient estimate is \( \widehat{\nabla J}_k(\boldsymbol{w}_{k, i-1}) = \nabla Q_k(\boldsymbol{w}_{k, i-1}; \boldsymbol{x}_{k, i}) \). Here, $\boldsymbol{\mu}_{k,i-1}$ is the step-size (or learning rate) of agent $k$ and $a_{k\ell}=[A]_{k,\ell}$ are the elements of the combination matrix. 

As in most gradient-based methods, the step-size in~\eqref{eq:atc} must be carefully tuned to balance the rate of convergence with the stability of the algorithm. Manual tuning (e.g., via a grid search) is costly in terms of time and computation, motivating the development of adaptive step-size strategies.

In centralized optimization, numerous adaptive methods have been proposed. As a representative example, consider Polyak's rule~\cite{polyak87} where the step-size is given by,
\begin{align}\label{eq:polyak}
    \mu_{k,i}=\frac{J(w_{k,i})-J(w^o)}{\lVert \nabla J_k(w_{k,i})\rVert^2}
\end{align}
Here, $J(w^o)$ is the minimum value of the cost function, which for many problems may be known in advance~\cite{boyd03}. Polyak's rule has been proven optimal for non-smooth convex optimization~\cite{polyak87}. Other centralized approaches include the Barzilai-Borwein method~\cite{barzilai88} and backtracking line search~\cite{armijo66}. Adaptive gradient methods, which rescale gradients coordinate-wise according to past gradient statistics, have also seen a surge in research interest due to their empirical success in machine learning. Included in this are Adam~\cite{kingma15}, AdaGrad~\cite{duchi11}, and their many variants.

Applying these centralized techniques to the decentralized setting is non-trivial. For example, an early decentralized implementation of Adam~\cite{nazari22} was shown in~\cite{chen22} to fail to converge to a stationary point under certain conditions due to step-size heterogeneity among agents. Other decentralized adaptive gradient methods have been implemented using techniques like dynamic consensus and variance reduction~\cite{chen22, huang24}.


Beyond adaptive gradient methods, decentralized counterparts of the Barzilai-Borwein rule have been studied in~\cite{gao21,iyanuoluwa22,hu21}. Similarly, the adaptive  step-size rule in~\cite{malitsky20} has been tailored for the decentralized setting in~\cite{zchen24,rikos25,ghaderyan25}. These approaches are not without their own challenges. The Barzilai-Borwein method has been shown to diverge for certain strictly-convex functions and has limited convergence guarantees~\cite{fletcher05}. Decentralized implementations of~\cite{malitsky20} appears unstable due to step-size heterogeneity~\cite{rikos25}. In~\cite{rikos25} this is tackled using a finite-time consensus protocol align both model parameters and agent step-sizes. While this ensures stability, it significantly increases communication cost.


A separate line of work considers backtracking line search in the decentralized setting~\cite{kuruzov24,kuruzov25,chen25}. This approach adaptively identifies the largest stable step-size, achieving linear convergence for strongly convex local objectives~\cite{kuruzov24,kuruzov25}. Importantly, it is parameter free, requiring no knowledge of smoothness or convexity constants. However, the repeated evaluation of local objective functions at each iteration makes line search computationally expensive and unsuitable when only stochastic gradients are available.

\subsection{Contribution}
Prior works incorporate step-size rules which have seen success in single-agent settings directly into decentralized strategies, by adjusting the local update step. In contrast, we derive an adaptive decentralized rule directly from first principles, accounting for the interaction between the optimal step-size choice and the graph topology. Our approach reduces to Polyak's rule in the single-agent case, hence we term our approach the graph-aware Polyak rule. The rule is parameter free, aside from requiring an estimate of the minimum objective value, analogously to the single-agent version of Polyak's rule. This value is frequently known or can be estimated a priori, for example in the overparamaterized learning setting. We further demonstrate that step-size heterogeneity is not necessarily harmful to decentralized algorithms.

\subsection{Notation}
Throughout the paper we use uppercase letters to denote matrices, lowercase for vectors and Greek symbols for scalars. Boldface symbols represent random quantities. For a matrix $Z$, the entry at row $k$ column $\ell$ is denoted $[Z]_{k,\ell}$. \textcolor{black}{The all-ones vector of size $K$ is denoted $\mathbbm{1}_K$ and the $M\times M$ identity matrix by $I_M$}. The Kronecker product is represented by $\otimes$, the Euclidean norm by $\norm{\cdot}$, and the trace operator by $\textrm{Tr}(\cdot)$.

\vspace{-5pt}\section{Algorithm Development}
We begin by rewriting the ATC-diffusion algorithm in~\eqref{eq:atc} more compactly in terms of network quantities:
\begin{align}\label{eq:atc-network}
    \w_i=\A(\w_{i-1}-\boldsymbol{U}_{i-1}\gradappcal(\w_{i-1}))
\end{align}
where ${\w_{i}\triangleq\textrm{col}\{\boldsymbol{w}_{k,i} \}}$, ${\gradappcal(\w_i)=\textrm{col}\{ \widehat{\nabla J}(\boldsymbol{w}_{k,i}) \}}$, the combination weights are defined by $\A\triangleq A\otimes I_M$ and the matrix of step-sizes $\boldsymbol{U}_{i-1}\triangleq\textrm{diag}\{\boldsymbol{\mu}_{k,i-1}\}\otimes I_M$. 

To continue, we make the following assumptions.
\begin{assumption}[Regularity conditions]\label{assump:regularity}
\textcolor{black}{The local loss functions $Q_k(\cdot;\cdot)$ are convex, so that for all \( y, z \in \mathbb{R}^M: \)}
\begin{align}
    Q_k(y;x_k)\geq Q_k(z;x_k)+\nabla Q_k(z;x_k)\transp(y-z) \label{eq:convex}
\end{align}
The local objective functions $J_k(\cdot)\triangleq \mathbb{E}Q_k(\cdot;\cdot)$ have Lipschitz gradients, which implies that:
\begin{align}\label{eq:assump-lipsch}
    J_k(x)\leq J_k(y) + \nabla J_k(y)\transp(x-y)+\delta_k\norm{x-y}^2
\end{align}
\end{assumption}
\textcolor{black}{We note that~\eqref{eq:convex} implies that the local objective functions $J_k(\cdot)$ are convex}.

\begin{assumption}[Combination Matrix]\label{assump:comb-mat}
The combination matrix $A$ is strongly-connected, symmetric, and doubly-stochastic.
\end{assumption}

To find the optimal choice of step-size at time \( i \), we will minimize the error obtained after a single iteration of the diffusion recursion. To this end, define the distance of the iterates to the optimal solution: 
\begin{align}
\phi(\w_{i})&\triangleq \norm{w^o\otimes \mathbbm{1}_K - \w_{i}}^2 \nonumber\\
&=\norm{\A(\werr_{i-1}+\boldsymbol{U}_{i-1}\gradappcal(\w_{i-1}))}^2 \label{eq:distance}
\end{align}
where $\werr_i\triangleq w^o\otimes \mathbbm{1}_K-\w_i$ is the error. 

Taking the derivative relative to the step-size: 
\begin{align} &\frac{\partial \phi(\w_{i+1})}{\partial \boldsymbol{\mu}_{k,i}} = \mathrm{Tr}\left[ \left(\nabla_{\boldsymbol{U}_i} \phi(\w_{i+1}\right)\transp\frac{\partial \boldsymbol{U}_i}{\partial\boldsymbol{\mu}_{k,i}}  \right]\nonumber\\
&= \mathrm{Tr}\left[ 2\A\A\big( \boldsymbol{U}_i\gradappcal(\w_i)+\werr_i\big)\gradappcal(\w_i)\transp \left(e_k e_k\transp\otimes I_M\right) \right] \label{eq:derivative}
\end{align}
where $e_k$ is the canonical basis vector with a one at position $k$, and zeros otherwise.

 For brevity we denote $B=AA$ (and $\mathcal{B}\triangleq \mathcal{A}\mathcal{A}$) to be the two-hop combination matrix, and its corresponding entries $b_{k\ell}\triangleq[B]_{k,\ell}$. We now proceed with simplifying the product in~\eqref{eq:derivative}, beginning with \textcolor{black}{the column vector}:
 \begin{align}
     &\left[\left(\boldsymbol{U}_i\gradappcal(\w_i)+\widetilde{\w}_i\right)\right]_{\ell} = \boldsymbol{\mu}_{\ell,i}\widehat{\nabla J}_\ell(\boldsymbol{w}_{\ell,i})+\boldsymbol{\widetilde{w}}_{\ell,i}
      \end{align}
      \textcolor{black}{Taking the outer product with the gradient}:
       \begin{align}
     &\left[ \left(\boldsymbol{U}_i\gradappcal(\w_i) + \widetilde{\w}_i\right) \gradappcal(\w_i) \transp\right]_{\ell,k} =\nonumber\\&\qquad\qquad \left(\boldsymbol{\mu}_{\ell,i}\widehat{\nabla J}_\ell(\boldsymbol{w}_{\ell,i})+\boldsymbol{\widetilde{w}}_{\ell,i}\right)\widehat{\nabla J}_k(\boldsymbol{w}_{k,i})\transp 
      \end{align}
      \textcolor{black}{Applying the two-hop combination matrix~$\B$, which averages over the columns in the preceding outer product}:
       \begin{align}
     &\left[ \B \left(\boldsymbol{U}_i\gradappcal(\w_i) + \widetilde{\w}_i\right) \gradappcal(\w_i) \transp\right]_{p,k}= \nonumber\\&~\qquad \sum_{\ell=1}^K b_{p\ell}(\boldsymbol{\mu}_{\ell,i}\widehat{\nabla J}_\ell(\boldsymbol{w}_{\ell,i}) + \boldsymbol{\widetilde{w}}_{\ell,i}) \widehat{\nabla J}_{k}(\boldsymbol{w}_{k,i})\transp \label{eq:matrix-subproducts}
 \end{align}
 
\setlength\arraycolsep{2.1pt}
\vspace{-5pt}\noindent The canonical basis vectors $e_k e_k\transp\otimes I_M$ in~\eqref{eq:derivative} act as column selection (or masking) matrices when post-multiplied to a matrix, setting all values in the matrix to zero except the  $k$-\textrm{th} block column. Thus,
 \begin{align} &\B\big(\boldsymbol{U}_i\gradappcal(\w_i)+\werr_i\big)\gradappcal(\w_i)\transp \left(e_k e_k\transp\otimes I_M\right)=\nonumber\\
  &\begin{bmatrix}  0 & \cdots & \sum_{\ell=1}^K b_{1\ell}(\boldsymbol{\mu}_{\ell,i}\widehat{\nabla J}_\ell(\boldsymbol{w}_{\ell,i})+\boldsymbol{\widetilde{w}}_{\ell,i})\widehat{\nabla J}_k(\boldsymbol{w}_{k,i})\transp & 0 & \cdots & \\ 
 0 & \cdots & \sum_{\ell=1}^K b_{2\ell}(\boldsymbol{\mu}_{\ell,i}\widehat{\nabla J}_\ell(\boldsymbol{w}_{\ell,i}) +\boldsymbol{\widetilde{w}}_{\ell,i})\widehat{\nabla J}_k(\boldsymbol{w}_{k,i})\transp & 0 &\cdots  \\
  &  & \vdots &  & 
\end{bmatrix} \end{align}\setlength\arraycolsep{5pt}
\textcolor{black}{As shown above, the only non-zero block column is the $k$-th $M\times M$ block. Applying the trace operator in~\eqref{eq:derivative} extracts the diagonal entries:}
\begin{align} &\frac{\partial \phi(\w_{i+1})}{\partial \boldsymbol{\mu}_{k,i}} 
\nonumber\\ &= \mathrm{Tr}\left( \sum_{\ell=1}^K b_{k\ell}(\boldsymbol{\mu}_{\ell,i}\widehat{\nabla J}_\ell(\boldsymbol{w}_{\ell,i}) +\boldsymbol{\widetilde{w}}_{\ell,i})\widehat{\nabla J}_k(\boldsymbol{w}_{k,i})\transp  \right)\nonumber\\
&=2\sum_{\ell=1}^Kb_{k\ell}\widehat{\nabla J}_k(\boldsymbol{w}_{k,i})\transp(\boldsymbol{\mu}_{\ell,i}\widehat{\nabla J}_\ell(\boldsymbol{w}_{\ell,i})+\boldsymbol{\widetilde{w}}_{\ell,i})
\end{align}
\textcolor{black}{where the last line follows from the identity ${\mathrm{Tr}(xy\transp)=y\transp x}$.}

Setting this to zero and solving for $\boldsymbol{\mu}_{k,i}$ we obtain a linear system of equations:
\begin{align}
    \sum_{\ell=1}^Kb_{1\ell}\widehat{\nabla J}_1(\boldsymbol{w}_{1,i})\transp\widehat{\nabla J}_\ell(\boldsymbol{w}_{\ell,i})\boldsymbol{\mu}_{\ell,i} &= \sum_{\ell=1}^Kb_{1\ell}\widehat{\nabla J}_1(\boldsymbol{w}_{1,i})\transp \boldsymbol{\widetilde{w}}_{\ell,i} \nonumber \\
    \sum_{\ell=1}^Kb_{2\ell}\widehat{\nabla J}_2(\boldsymbol{w}_{2,i})\transp\widehat{\nabla J}_\ell(\boldsymbol{w}_{\ell,i})\boldsymbol{\mu}_{\ell,i} &= \sum_{\ell=1}^Kb_{2\ell}\widehat{\nabla J}_2(\boldsymbol{w}_{2,i})\transp \boldsymbol{\widetilde{w}}_{\ell,i} \nonumber \\
    &\vdots \nonumber
\end{align}
Or, in matrix form:
\setlength\arraycolsep{2.1pt}
\begin{subequations}
\begin{align} &\hspace{0.4\linewidth} \boldsymbol{Z}_i \boldsymbol{u}_i = - \boldsymbol{c}_i \label{eq:lse} \\
&\textrm{where:} \nonumber\\
    &\boldsymbol{Z}_i\!=\!\small\!\begin{bmatrix}b_{11}\widehat{\nabla J}_1(\boldsymbol{w}_{1,i})\transp\widehat{\nabla J}_1(\boldsymbol{w}_{1,i}) &b_{12}\widehat{\nabla J}_1(\boldsymbol{w}_{1,i})\transp\widehat{\nabla J}_2(\boldsymbol{w}_{2,i}) &\cdots \\ b_{21}\widehat{\nabla J}_2(\boldsymbol{w}_{2,i})\transp\widehat{\nabla J}_1(\boldsymbol{w}_{1,i}) &b_{22}\widehat{\nabla J}_2(\boldsymbol{w}_{2,i})\transp\widehat{\nabla J}_2(\boldsymbol{w}_{2,i}) &\cdots\\& \vdots & \end{bmatrix}\normalsize \\
    & \hspace{0.435\linewidth} \boldsymbol{u}_i = \begin{bmatrix}\boldsymbol{\mu}_{1,i}\\ \boldsymbol{\mu}_{2,i} \\\vdots\end{bmatrix} \\
    & \hspace{0.28\linewidth} \boldsymbol{c}_i = \begin{bmatrix} \sum_{\ell}b_{1\ell}\widehat{\nabla J}_1(\boldsymbol{w}_{1,i})\transp\boldsymbol{\widetilde{w}}_{\ell,i} \\ \sum_{\ell}b_{2\ell}\widehat{\nabla J}_2(\boldsymbol{w}_{2,i})\transp\boldsymbol{\widetilde{w}}_{\ell,i} \\ \vdots\end{bmatrix}
\end{align}\end{subequations}
The optimal step-size choice hence reduces to the solution of a decentralized linear system of equations. In other words, each agent $k$ aims to solve the $k$-\textrm{th} row of the system of equations which requires quantities from agents within its two-hop neighborhood ($b_{k\ell}=0$ if $\ell$ is not a neighbor of agent $k$'s neighbors). To solve~\eqref{eq:lse} we use an iterative method based on the randomized block Kaczmarz method~\cite{elfving80}:
\begin{subequations}\label{eq:kaczmarz}
\begin{align} &\boldsymbol{\mu}^{(t)}_{k,i}
= \boldsymbol{\mu}_{k,i}^{(t-1)}-\boldsymbol{\alpha}_k \sum_{\ell=1}^Kb_{k\ell}\widehat{\nabla J}_k(\boldsymbol{w}_{k,i})\transp\boldsymbol{\widetilde{w}}_{\ell,i} \nonumber\\&\quad-\boldsymbol{\alpha}_k\sum_{\ell=1}^K b_{k\ell}\boldsymbol{\mu}_{\ell,i}^{(t-1)}\widehat{\nabla J}_k(\boldsymbol{w}_{k,i})\transp\widehat{\nabla J}_\ell(\boldsymbol{w}_{\ell,i}) \\
&\textrm{where,} \nonumber\\
&\boldsymbol{\alpha}_k = \frac{b_{kk}\lVert \widehat{\nabla J}_k(\boldsymbol{w}_{k,i})\rVert^2}{\sum_{\ell=1}^K \left( b_{k\ell}\widehat{\nabla J}_k(\boldsymbol{w}_{k,i})\transp\widehat{\nabla J}_\ell(\boldsymbol{w}_{\ell,i}) \right)^2} \label{eq:alpha_k}
\end{align}\end{subequations}
initialized with $\boldsymbol{\mu}_{k,i}^{(0)}=0$. 

The iterative algorithm for finding $\boldsymbol{\mu}_{k,i}$ in~\eqref{eq:kaczmarz} requires knowledge of the optimal solution, ${\boldsymbol{\widetilde{w}}_{\ell,i}\triangleq w^o-\boldsymbol{w}_{\ell,i}}$. To avoid this, we replace~$\boldsymbol{\mu}_{k,i}$ by the upper bound:
\begin{align} &\widehat{\nabla J}_k(\boldsymbol{w}_{k,i}) \transp\boldsymbol{\widetilde{w}}_{\ell,i} = \widehat{\nabla J}_k(\boldsymbol{w}_{k,i})\transp(w^o-\boldsymbol{w}_{\ell,i}) \nonumber\\
&= \widehat{\nabla J}_k(\boldsymbol{w}_{k,i})\transp(w^o-\boldsymbol{w}_{k,i}+\boldsymbol{w}_{k,i}-\boldsymbol{w}_{\ell,i}) \nonumber\\
&= \widehat{\nabla J}_k(\boldsymbol{w}_{k,i})\transp(w^o-\boldsymbol{w}_{k,i})+\widehat{\nabla J}_k(\boldsymbol{w}_{k,i})\transp(\boldsymbol{w}_{k,i}-\boldsymbol{w}_{\ell,i}) \nonumber\\
&\stackrel{(a)}{\leq} \widehat{J}_k(w^o) - \widehat{J}_k(\boldsymbol{w}_{k,i})+\widehat{\nabla J}_k(\boldsymbol{w}_{k,i})\transp(\boldsymbol{w}_{k,i}-\boldsymbol{w}_{\ell,i}) \nonumber\\
&\stackrel{(b)}{\leq} \widehat{J}_k(w^o) - \widehat{J}_k(\boldsymbol{w}_{k,i})+ \widehat{J}_k(\boldsymbol{w}_{k,i})-\widehat{J}_k(\boldsymbol{w}_{\ell,i})\nonumber\\
&\qquad \qquad+\delta_k\lVert \boldsymbol{w}_{k,i}-\boldsymbol{w}_{\ell,i}\rVert^2 \nonumber\\
&= \widehat{J}_k(w^o) -\widehat{J}_k(\boldsymbol{w}_{\ell,i})+\delta_k\lVert \boldsymbol{w}_{k,i}-\boldsymbol{w}_{\ell,i}\rVert^2 \label{eq:upperbound}
\end{align}
where $(a)$ and $(b)$ follow from Assumption~\ref{assump:regularity} by convexity and Lipschitz smoothness, respectively. Note that $J_k(\boldsymbol{w}_{\ell,i})$ corresponds to the stochastic cost of agent \( k \), evaluated with the parameter vector of agent \( \ell \).

Equation~\eqref{eq:upperbound} provides an upper bound for $\boldsymbol{\mu}_{k,i}$. To ensure the step-size remains below this bound, we replace it by the more conservative bound:
\begin{align}
\widehat{\nabla J}_k(\boldsymbol{w}_{k,i}) \transp\boldsymbol{\widetilde{w}}_{\ell,i}\lesssim \widehat{J}_k(w^o) - \widehat{J}_k(\boldsymbol{w}_{\ell,i}) 
\end{align}
\textcolor{black}{The step-size update becomes}:
\begin{align}
    \boldsymbol{\mu}_{k,i}^{(t)} &= \boldsymbol{\mu}_{k,i}^{(t-1)} - \boldsymbol{\alpha}_k \sum_{\ell=1}^K b_{k\ell}(\widehat{J}_k(w^o)-\widehat{J}_k(\boldsymbol{w}_{\ell,i})) \nonumber \\ &\qquad- \boldsymbol{\alpha}_k\sum_{\ell=1}^K b_{k\ell}\boldsymbol{\mu}_{\ell,i}^{(t-1)}\widehat{\nabla J}_k(\boldsymbol{w}_{k,i})\transp\widehat{\nabla J}_\ell(\boldsymbol{w}_{\ell,i}) \label{eq:step-size-rule}
\end{align}
After $t=T$ iterations we set $\boldsymbol{\mu}_{k,i}=\boldsymbol{\mu}_{k,i}^{(T)}$. In practice we find that one to two iterations of the Kaczmarz method suffice for calculating the step-size with sufficient accuracy. 

The ATC-diffusion algorithm with our graph-aware Polyak step-size is summarized in Algorithm~\ref{alg:gap}.  \textcolor{black}{Where possible we have split the update from~\eqref{eq:step-size-rule} into intermediate quantities using the one-hop combination matrix and also include a non-negativity constraint}.

\begin{algorithm}
\caption{ATC-Diffusion with Graph-Aware Polyak's Rule}\label{alg:gap}
\begin{algorithmic}
\For{\texttt{i=1,2,\dots}}
    \State $g_{k,i} \gets w_{k,i-1} - \mu_{k,i-1}\widehat{\nabla J}_k(w_{k,i-1})$
    \State $w_{k,i} \gets \sum_{j\in\mathcal{N}_k} a_{kj}g_{j,i}$
    \State $\mu_{k,i}^{(0)}\gets 0$
    \State $\alpha_k $ \texttt{ from \eqref{eq:alpha_k}} 
    \State $\beta_k \gets \alpha_k\sum_{j=1}^K b_{kj}(\widehat{J}_k(w^o) - \widehat{J}_k(w_{j,i}))$
    \For{\texttt{t=1,\dots,T}}
            \State $v_{k,i}^{(t)} \gets \sum_{\ell=1}^K a_{k\ell}\mu_{\ell,i}^{(t-1)}\widehat{\nabla J}_\ell(w_{\ell,i})$
            \State $ \mu^{(t)}_{k,i} \gets \mu_{k,i}^{(t-1)}-\beta_k - \alpha_k\widehat{\nabla J}_k(w_{k,i})\transp\sum_{\ell=1}^K a_{k\ell}v_{\ell,i}^{(t)} $
        \EndFor
        \State $\mu_{k,i}\gets \textrm{max}\left(\mu_{k,i}^{(T)},0\right)$
      \EndFor
\end{algorithmic}
\end{algorithm}

As with Polyak's rule, the step-size in~\eqref{eq:step-size-rule} requires evaluating the local objective $J_k(\cdot)$ at the optimal solution $w^o$. Here, however, $w^o$ is the global minimizer across the network, rather than the local minimizer of agent $k$'s objective. When agents have homogeneous data the global and local minimizers coincide, $(w_k^o=w^o)$ and for many problems, the local minimum $J_k(w_k^o)$ may be known \textit{a priori}~\cite{boyd03}. In heterogeneous settings, the cost at the global minimum may still be known in advance, for example, in over-parameterized models where $J_k(w^o)=0$.

Finally, we demonstrate the limiting behavior of~\eqref{eq:step-size-rule} in the single-agent case ($K=1$) case, assuming deterministic evaluations of the objective and gradient. We have coefficient $\alpha_k=1$ and the update becomes:
\begin{align}\mu_{i}^{(t)}&= \mu_{i}^{(t-1)}-(J(w^o)-J(w_{i}))-\mu_{i}^{(t-1)}\nabla J(w_{i})\transp\nabla J(w_{i})
\end{align}
Taking the limit as $t\rightarrow\infty$:
\begin{align}\mu_{i}^{(\infty)}&= \mu_{i}^{(\infty)}-(J(w^o)-J(w_{i}))-\mu_{i}^{(\infty)}\nabla J(w_{i})\transp\nabla J(w_{i})
\end{align}
Solving for $\mu_{i}^{(\infty)}$ yields $\mu_{i}^{(\infty)}=\frac{J(w_i)-J(w^o)}{\norm{\nabla J(w_i)}^2}$ which corresponds to Polyak's rule.\hfill $\square$\\

\section{Numerical Results}
\begin{figure}[htb]
  \centering
  \centerline{\includegraphics[width=7.8cm]{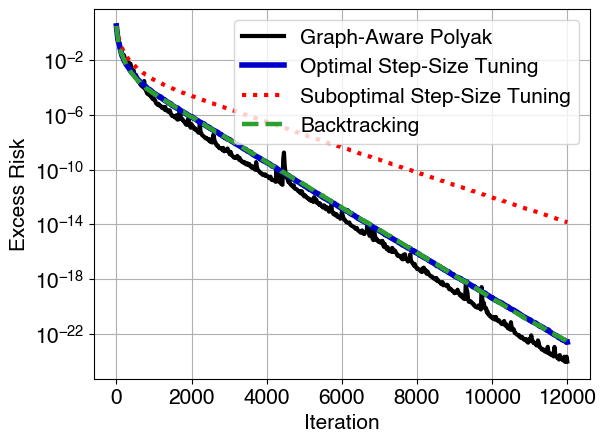}}
\caption{\textcolor{black}{Excess-Risk for a linear regression problem on an Erdős–Rényi Graph}}
\label{fig:results}
\end{figure}

We demonstrate our step-size rule for a linear regression problem. Each agent holds samples generated via a linear model, $\gamma_{k,n}=h_{k,n}\transp w^o+\nu_{k,n}$  and attempts to fit a model using a least squares objective,
\begin{align}
    J_k(w)=\frac{1}{2N}\sum_{n=1}^N(\gamma_{k,n}-h_{k,n}\transp w)^2
\end{align}
The features, $h_{k,n}$, and noise, $\nu_{k,n}$, follow normal distributions. We demonstrate our method on an Erdős–Rényi graph with 8 agents in the over-parameterized setting. Each agent employs stochastic evaluations of the objective by selecting a random sample $n_i$ and computing both the objective value $\widehat{J}(\cdot)$ and its gradient $\widehat{\nabla J}_k(\cdot)$. We set the number of iterations in the inner loop in Algorithm~\ref{alg:gap} to $T=2$. Figure~\ref{fig:results} compares the proposed step-size rule with ATC-Diffusion run with the best fine-tuned step-size. For reference, we also include a suboptimal step size set to half the fine-tuned one as well as the \textcolor{black}{We have also adapted the backtracking rule from~\cite{kuruzov25} to ATC-diffusion to enable even comparison}. The adaptive and fine-tuned algorithms exhibit similar performance, demonstrating the effectiveness of the proposed rule without the need for fine-tuning. \textcolor{black}{Although the graph-aware Polyak rule requires additional computation and communication, it obviates step-size tuning, avoiding stability issues and slow convergence.} \textcolor{black}{We highlight two key observations in the work, namely that the optimal local step-size in decentralized algorithms is influenced by the network topology, and that step-size heterogeneity is not always detrimental to performance}.

\vspace{-15pt}
\section{Conclusion}

We have derived an adaptive step-size method for decentralized optimization directly from first principles. Compared to prior art, the resulting step-size rule incorporates the effect of network topology in the design. The method is parameter free apart from requiring knowledge of the cost at the global minimizer (which may be known in advance), and allows for heterogeneous step-sizes across the network. Numerical experiments demonstrate that the method performs comparably to fine-tuning the step-size.

\section{References}
\renewcommand{\refname}{}
\setlength{\bibsep}{0pt plus 0.3ex}
\bibliographystyle{IEEEbib}
\vspace{-20pt}
\bibliography{References}

\end{document}